\documentclass[conference]{IEEEtran}
\IEEEoverridecommandlockouts
\usepackage{mathpazo}
\usepackage{times}
\usepackage{amsmath}
\usepackage{amsfonts}
\usepackage{latexsym}
\usepackage{amssymb}
\usepackage{mathabx}
\usepackage{float}
\usepackage{wrapfig}
\usepackage{upref}
\usepackage{theorem}
\usepackage{graphicx}
\usepackage{psfrag}
\usepackage{cite}
\usepackage{textcomp}
\usepackage{xcolor}

\usepackage{color}
\usepackage{tikz}
\usetikzlibrary{arrows.meta,bending}
\makeatletter
\newcommand\curvearrowed@[3]
  {%
    \begin{tikzpicture}
      \node[inner sep=.05ex](a){\kern-.25ex#3\kern-.05ex};
      \draw[arrows={-Latex[#1]},line width=#2]
        (a.north west) to[bend left=45] (a.north east);
    \end{tikzpicture}%
  }
\newcommand\curvearrowed[1]
  {%
    \relax
    \ifmmode
      \mathchoice
        {\curvearrowed@{}{.4pt}{$\displaystyle #1$}}
        {\curvearrowed@{}{.4pt}{$\textstyle #1$}}
        {\curvearrowed@{scale=.8}{.325pt}{$\scriptstyle #1$}}
        {\curvearrowed@{scale=.7}{.25pt}{$\scriptscriptstyle #1$}}%
    \else
      \curvearrowed@{}{.4pt}{#1}%
    \fi
  }
\makeatother

\usepackage{algorithm,algpseudocode}

\makeatletter
\newcommand{\removelatexerror}{\let\@latex@error\@gobble}
\makeatletter
\newcommand{\proofpart}[2]{%
	\par
	\addvspace{\medskipamount}%
	\noindent\emph{Part #1: #2}\par\nobreak
	\addvspace{\smallskipamount}%
	\@afterheading
}
\makeatother

\theoremstyle{plain}
\theorembodyfont{\normalfont\slshape}

\newtheorem{thm}{Theorem$\!$}
\newenvironment{theorem}
{\begin{thm}\hspace*{-1ex}{\bf.}}{\end{thm}}

\newtheorem{clm}[thm]{Claim$\!$}
\newenvironment{claim}{\begin{clm}\hspace*{-1ex}{\bf.}}{\end{clm}}

\newtheorem{lem}[thm]{Lemma$\!$}
\newenvironment{lemma}{\begin{lem}\hspace*{-1ex}{\bf.}}{\end{lem}}

\newtheorem{prop}[thm]{Proposition$\!$}

\newtheorem{cor}[thm]{Corollary$\!$}

\newtheorem{defn}[thm]{Definition$\!$}
\newenvironment{definition}{\begin{defn}\hspace*{-1ex}{\bf.}}{\end{defn}}

\newtheorem{xmpl}[thm]{Example$\!$}
\newenvironment{example}{\begin{xmpl}\hspace*{-1ex}{\bf.}}{\hfill $\Box$ \end{xmpl}}

\newtheorem{cnstr}{Construction$\!$}

\newtheorem{rmk}[thm]{Remark$\!$}

\newcounter{enumrom}
\renewcommand{\theenumrom}{(\roman{enumrom})}

\makeatletter
\renewcommand{\@endtheorem}{\endtrivlist}
\makeatother

\makeatletter
\renewcommand{\thefigure}{{\@arabic\c@figure}}
\renewcommand{\fnum@figure}{{\bf Figure\,\thefigure}}
\makeatother

\renewcommand{\le}{\leqslant}
\renewcommand{\leq}{\leqslant}
\renewcommand{\ge}{\geqslant}
\renewcommand{\geq}{\geqslant}

\newcommand{\Cref}[1]{Co\-ro\-lla\-ry\,\ref{#1}}

\outer\def\proclaim #1. #2\par{\medbreak
 \noindent{\bf#1.\enspace}{\sl#2\par}%
 \ifdim\lastskip<\medskipamount \removelastskip\penalty55\medskip\fi}

\begin{document}

\title{Balanced and Swap-Robust Trades for Dynamical Distributed Storage
\thanks{This work was supported by the NSF grants CCF 1816913 (CJC) and 1814298 (OM). The authors acknowledge useful discussions with Dylan Lusi.}}

\makeatletter
\newcommand{\linebreakand}{%
  \end{@IEEEauthorhalign}
  \hfill\mbox{}\par
  \mbox{}\hfill\begin{@IEEEauthorhalign}
}
\makeatother

\author{\IEEEauthorblockN{Chao Pan}
\IEEEauthorblockA{\textit{ECE Department} \\
\textit{University of Illinois}\\
chaopan2@illinois.edu} 
\and
\IEEEauthorblockN{Ryan Gabrys}
\IEEEauthorblockA{\textit{Calit2} \\
\textit{University of California-San Diego} \\
ryan.gabrys@gmail.com}
\and
\IEEEauthorblockN{Xujun Liu}
\IEEEauthorblockA{\textit{Foundational Mathematics Department} \\
\textit{Xi'an Jiaotong-Liverpool University}\\
xujun.liu@xjtlu.edu.cn}
\linebreakand
\IEEEauthorblockN{Charles Colbourn}
\IEEEauthorblockA{\textit{School of CAI} \\
\textit{Arizona State University}\\
charles.colbourn@asu.edu}
\and
\IEEEauthorblockN{Olgica Milenkovic}
\IEEEauthorblockA{\textit{ECE Department} \\
\textit{University of Illinois}\\
milenkov@illinois.edu}
}
\maketitle

\begin{abstract} Trades, introduced by Hedayat~\cite{hedayat1990theory}, are two sets of blocks of elements which may be exchanged (traded) without altering the counts of certain subcollections of elements within their constituent blocks. They are of importance in applications where certain combinations of elements dynamically become prohibited from being placed in the same group of elements, since in this case one can trade the offending blocks with allowed ones. This is particularly the case in distributed storage systems, where due to privacy and other constraints, data of some groups of users cannot be stored together on the same server. We introduce a new class of \emph{balanced trades}, important for access balancing of servers, and \emph{perturbation resilient balanced trades}, important for studying the stability of server access frequencies with respect to changes in data popularity. The constructions and bounds on our new trade schemes rely on specialized selections of defining sets in minimal trades and number-theoretic analyses.
\end{abstract}


\section{Introduction}

An important problem in distributed storage is to balance access requests to servers in 
order to prevent service time bottlenecks. Access balancing aims to limit the number of server session logon requests and queues excess requests until server resources become available. As data access directly relates to the popularity of items on a server, it is important to design storage systems that balance the popularity of data chunks on each server. 

Existing distributed storage systems categorize data into hot, warm and cold based on their popularity and use different 
mechanisms for storing them~\cite{cherkasova2004analysis,rawat2016locality,joshi2014delay}. High-level partitions like these fail to account for popularity-variability within each class and the resulting dynamic changes in data demands, and they do not couple the balancing process with distributed storage repair solutions~\cite{dimakis2010network,el2010fractional,silberstein2015optimal}. As a result, although issues such as delay-storage tradeoffs, volume (load) balancing, and chunk allocation have been studied in depth~\cite{leong2012distributed,joshi2014delay,aktacs2021evaluating}, access balancing has received significantly less attention. To address this problem, the authors of~\cite{dau2018maxminsum} introduced \emph{MaxMin Steiner triple systems} and related combinatorial designs for fractional repetition distributed storage coding~\cite{el2010fractional}. MaxMin Steiner triple systems maximize the minimum average popularity score of items stored on different servers in the system and thereby ensure control of the discrepancy of average server popularity scores~\cite{brummond2019kirkman,chee2020access,colbourn2021egalitarian}. These systems are constructed by selecting a specific packing or design, and associating data items of given popularities with elements of the packing or design. The MaxMin Steiner balancing approach can also be integrated with batch codes~\cite{el2010fractional,silberstein2015optimal,silberstein2016optimal}.

In~\cite{hedayat1990theory}, the following question was posed: ``In the application of $t$-designs we may be confronted with a situation where some blocks become too costly to be selected for experimentation. For example, combining certain tasks within one block may be unacceptable to the experimenter. Suppose the available $t$-designs in the literature contain such blocks and by renaming we cannot dispose of these undesirable blocks. The theory of trade-offs can tell us how to trade undesirable blocks with those which are acceptable to the experimenter.''  In access balancing such situations arises when one needs to redistribute files or data chunks on servers to account for constraints on items which cannot be stored together or to mitigate the effect of dynamical changes in the popularity of items.

This work demonstrates that such constrained trading problems may be resolved through the use of new classes of balanced and swap-robust trades. A trade with parameters $(v,k,t)$ is a pair of disjoint sets of blocks $T^{(1)},T^{(2)}$ of equal size over a ground set of $v$ points such that all blocks are of size $k$ and any subset of $t$ points appears the same number of times in blocks of $T^{(1)}$ and blocks of $T^{(2)}$. The set of blocks in $T^{(1)}$ and $T^{(2)}$ are subgroups of blocks in combinatorial designs which are used for distributed and batch code constructions. Replacing blocks in $T^{(1)}$ by blocks in $T^{(2)}$ when the need arises does not change the global block intersection properties. In the setting of access-balanced storage with dynamically changing item popularities, undesirable blocks may represent blocks whose average popularity has decreased or increased significantly or blocks who are dynamically faced with storing disallowed combinations of data items. Importantly, for dynamic access balancing, blocks in $T^{(1)}$ can be traded by blocks in $T^{(2)}$ without changing the properties of the underlying design but changing the points (e.g., data chunks or files) allocated to the blocks (e.g., servers). Due to space limitations, we relegate the detailed description of how to use trades in distributed storage system to the journal version of this work. 

The first question of interest is: Can trades themselves be made ``balanced''? Labeling the points with numbers in $\{{1,2,\ldots,v\}}$, the question of balancing translates to ensuring that all blocks in $T^{(1)}$ and $T^{(2)}$ have the same sum of label values. Our first result is that minimal trades (i.e., trades with the smallest number of blocks) can be easily balanced (see Section~\ref{sec:balanced-trades}). Balancing is needed even when no changes in popularities are present to ensure maintaining near-uniform server access; minimal trades are of interest since they require replacing (trading) the smallest number of blocks to achieve the desired rebalancing or to meet a given placement constraint. The second question is how to design balanced trades that are as robust as possible to limited-magnitude swaps in popularities of items? For $v=4(t+1)$, where $t$ is a positive integer, we show that there exist balanced trades whose popularity sums do not change more than roughly $2(t+1)$ when all popularity values are allowed to undergo arbitrary \emph{adjacent} swaps. Interestingly, for $v=16$ and $v=20$ the corresponding balanced minimal trades are unique.  

\section{Basics of Trades} \label{sec:trades}

We start our exposition with some important definitions.

Given a set of elements $V=[v] = \{1,\ldots,v\}$, let $P(v)$ denote the power set of $V$ and let $P_k(v) \subseteq P(v)$ be the set of all subsets (blocks) of $P(v)$ that have cardinality $k,$ for some positive integer $k < v$. 

\begin{definition}
A $(v,k,t)$ trade is a pair of sets $\{{T^{(1)}, T^{(2)}\}}$,  $T^{(i)} \subseteq P_k(v), i=1,2,$ such that $|T^{(1)}| = |T^{(2)}|$ and 
$T^{(1)} \cap T^{(2)} = \emptyset$, satisfying the following property: For any $B_t \in P_t(v)$, the number of blocks in $T^{(1)}$ that contain $B_t$ is the same as the number of blocks in $T^{(2)}$ that contain $B_t$.
The volume of the trade is the number of blocks in $T^{(1)}$ and $|T^{(1)}| = |T^{(2)}|$. A $(v,k,t)$ trade is 
termed \emph{minimal} if it has the smallest possible number of blocks. 
\end{definition}

We start by describing a construction of minimal trades based on Hwang~\cite{hwang1986structure} which allows one to balance the sum of elements in blocks (i.e., ensure discrepancy = 0) of both $T^{(1)}$ and $T^{(2)},$ using specialized choices of the sets defining the blocks of the trades. We describe the construction and an accompanying new proof that immediately reveals how to perform balancing (note that no proof for this result was provided in~\cite{hwang1986structure}).
For each $\ell \in \{1, 2, \ldots, n\}$, let 
$$P_{\ell} = \{(i_1, i_1 + 1) \cdots (i_{\ell}, i_{\ell} + 1) \, | \, \{i_1, \ldots, i_{\ell}\} \subseteq \{1,\ldots, 2n-1\}\},$$
be a set of permutations over a set of $2n$ elements, where $(i_1, i_1 + 1)(i_2, i_2 + 1) \ldots (i_{\ell}, i_{\ell} + 1)$ denotes the permutation which swaps $i_j$ with $i_j + 1$, for $j = 1, \ldots, \ell,$ and leaves all other elements fixed. Furthermore, let $P_0 = \{(e)\}$ consist of the identity permutation $e$. Also, define 
$$\Delta_{2n} = \bigcup\limits_{\ell \text{ even}} P_{\ell} \quad \text{ and } \quad \bar{\Delta}_{2n} = \bigcup\limits_{\ell \text{ odd}} P_{\ell}.$$

\begin{theorem}(\cite{hwang1986structure})\label{thm:hwang}
There exists a $(v,k,t)$ trade of volume $2^t$. This volume is the smallest (minimal) possible volume for the given choice of parameters.
\end{theorem}

\begin{IEEEproof}
We provide a simple proof for the first claim. Let $S_1, S_2, \ldots, S_{2t + 3}$ be subsets of the point set $V$ that can be partitioned into $t+1$ pairs of the form 
$$S_1, S_2; S_3, S_4; \ldots S_{2t+1}, S_{2t+2},$$ with the addition of one unpaired set $S_{2t+3}$. We henceforth refer to the sets as \emph{defining sets} and the pairs as \emph{companions}. Assume that the defining sets have the following properties: 
\begin{enumerate}
\item $S_i \cap S_j = \emptyset$ for $i \neq j$;
\item $|S_{2i-1}| = |S_{2i}| \ge 1$, for $i = 1, \ldots, t+1$;
\item $\sum\limits_{i = 1}^{t+2} |S_{2i-1}| = k$.
\end{enumerate}
Respectively, define $T = \{T^{(1)}, T^{(2)}\}$ as
$$\{{ (S_{\sigma(1)} \cup S_{\sigma(3)} \cup \ldots \cup S_{\sigma(2t+1)} \cup S_{2t+3}): \sigma \in \Delta_{2t+2} \}},$$
$$\{{(S_{\bar{\sigma}(1)} \cup S_{\bar{\sigma}(3)} \cup \ldots \cup S_{\bar{\sigma}(2t+1)} \cup S_{2t+3}): \bar{\sigma} \in \bar{\Delta}_{2t+2}\}}.$$
We used $(S_{\sigma(1)} \cup S_{\sigma(3)} \cup \ldots \cup S_{\sigma(2t+1)} \cup S_{2t+3})$ and $(S_{\bar{\sigma}(1)} \cup S_{\bar{\sigma}(3)} \cup \ldots \cup S_{\bar{\sigma}(2t+1)} \cup S_{2t+3})$ to denote blocks whose elements represent the unions of the corresponding sets $S$.

\begin{claim}
The trade $T$ has volume $2^t$.
\end{claim}
\begin{IEEEproof}
Since $\Delta_{2t+2}$ represents the set of all permutations with an even number of transpositions, for $t$ even we have 
$$|T^{(1)}| = |\Delta_{2t+2}| = {t+1 \choose 0} + {t+1 \choose 2} + \ldots + {t+1 \choose t} = 2^t,$$
$$|T^{(2)}| = |\bar{\Delta}_{2t+2}| = {t+1 \choose 1} + {t+1 \choose 3} + \ldots + {t+1 \choose t+1} = 2^t.$$
The same results hold for odd $t$.
\end{IEEEproof}

\begin{claim}
$T$ is a $(v,k,t)$-trade.
\end{claim}

\begin{IEEEproof}
We show that every block has size $k$ and that every $t$-subset of $V$ is contained in the same number of blocks in $T^{(1)}$ and $T^{(2)}$.

Firstly, the claim that every block has size $k$ follows from Property (3) of the collection of sets $S_1, S_2, \ldots, S_{2t + 3}$. 
Secondly, to show that every $t$-subset of $V$ is contained in the same number of blocks in $T^{(1)}$ and $T^{(2)}$, let $U$ be a $t$-subset of the ground set $V$.

\textbf{Case 1:} There exists some $i \in \{1, \ldots, t+1\}$ such that $U \cap S_{2i-1} \neq \emptyset$ and $U \cap S_{2i} \neq \emptyset$. Then, the construction of the trade ensures that $U$ is not contained in any block.

\textbf{Case 2:} For each pair $S_{2i-1}, S_{2i}$, $1 \le i \le t+1$, $U$ has a nonempty intersection with at most one of the sets. Without loss of generality, assume that $U \cap S_{2 i - 1} \neq \emptyset$ for indices $i \in \{{i_1, i_2, \ldots, i_{h_1}\}}$ and $U \cap S_{2 j} \neq \emptyset$ for indices $j \in \{{k_1, k_2, \ldots, k_{h_2}\}}$. Furthermore, without loss of generality, assume that $(h_1+h_2)$ is odd and that $t$ is even. Then $U$ is contained in 
$${t+1-h_1-h_2 \choose  1} + \ldots + {t+1-h_1-h_2 \choose  t-h_1-h_2}$$ 
blocks of $T^{(1)}$ and  
$${t+1-h_1-h_2 \choose  0} + \ldots + {t+1-h_1-h_2 \choose  t+1-h_1-h_2}$$
blocks of $T^{(2)}$. In both cases, the sums equal $2^{t-h_1-h_2}$. This establishes the claim.
\end{IEEEproof}
This completes the proof of the theorem.
\end{IEEEproof}
We remark that all trades of volume $2^t$ must have a structure as described in the above proof, as shown in~\cite{hwang1986structure}.

\section{Balanced Minimal Trades} \label{sec:balanced-trades}

We introduce next the notion of block-sums and block-discrepancies needed for our subsequent discussion.
\begin{definition}
Let $B=(b_1,b_2,\ldots,b_k)$ be a block in $T^{(i)}$, where $i \in \{{1,2\}}$. The block-sum of $B$ equals $\Sigma_B=\sum_{i=1}^k\,b_i$, while the minimum and maximum block sums of $T^{(i)}$ are defined as $\min_{B \in T^{(i)}}\, \Sigma_B$ and 
$\max_{B \in T^{(i)}}\, \Sigma_B$, for $i \in \{{1,2\}}$. The block-discrepancy of $T^{(i)}$ is defined as 
$\max_{B \in T^{(i)}}\, \Sigma_B - \min_{B \in T^{(i)}}\, \Sigma_B$, for $i\in\{{1,2\}}$.
\end{definition}
Henceforth, we tacitly assume that the points are labeled according to their popularity score so that no two popularities are the same: Label ``$1$'' indicates the most popular data chunk, while label ``$v$'' denotes the least popular data chunk.

We ask the following questions: Can we ensure that both $T^{(1)}$ and $T^{(2)}$ have zero block-discrepancy? We affirmatively answer this question in the following lemma.

\begin{lemma} \label{lem:baltrade}
There exists a minimal $(v,k,t)$ trade $T = \{T^{(1)}, T^{(2)}\}$ of volume $2^t$ with $v=4(t+1)$ points such that both $T^{(1)}$ and $T^{(2)}$ are perfectly balanced: Every block has the same size $k=2(t+1)$ and block-sum equal to $(t+1)(4t+5)$. 
\end{lemma}

\begin{IEEEproof}
Let $T = \{T^{(1)}, T^{(2)}\}$ be defined as in the previous theorem, and select the companion sets as
$$S_1 = \{1,4\}, S_2 = \{2,3\}, S_3 = \{5,8\}, S_4 = \{6,7\}, \ldots, $$
$$S_{2t+1} = \{4t+1,4t+4\}, S_{2t+2} = \{4t+2,4t+3\},\; S_{2t+3}=\emptyset.$$
Clearly, this choice of sets satisfies the conditions required of the construction to result in a trade. It is also easy to see that every block has the same size $2t+2$ and that the sum of the labels in each block equals $(t+1) (4t+5)$.
\end{IEEEproof}
Another obvious choice for the defining sets $S_1,S_2,\ldots, S_{2t+2}$ for doubly-even $v$ is
$$S_1 = \{1,v-1\}, S_2 = \{2,v-2\}; S_3 = \{3,v-3\}, S_4 = \{4,v-4\};$$
$$\ldots S_{2t+1} = \{2t+1,v-2t-1\}, S_{2t+2} = \{2t+2, v-2t-2\}.$$ 
We determine next how many \emph{balanced} partitions into pairs of companion sets of cardinality two, i.e. partitions such that $\sum_{\ell \in S_{2i-1} } \ell=\sum_{\ell \in S_{2i} }\ell$, where $i=1,2,\ldots,t+1$. To the best of the authors' knowledge, this question has not been previously addressed in the literature and appears to be difficult. Computer search shows that for $v=12,16,20,24$ the number of valid partitions equals $86,1990,74323,4226026$. A straightforward, yet loose, lower bound, equals $\frac{(2(t+1))!}{(t+1)!\,2^{t+1}}$. This follows as one can first choose a pair of integers in $\{{1,2,\ldots,2(t+1)\}}$, say $a$ and $b$, and then select $a+2(t+1)$ and $b+2(t+1)$ in $\{{2(t+1)+1,2(t+1)+2,\ldots,4(t+1)\}}$ as these four elements constitute two balanced companion sets. 

An upper bound can be obtained by viewing the problem of counting balanced set partitions as a problem pertaining to integer partitions. Specifically, the sum of all elements in $\{{1,2,\ldots,v\}}$ equals $s=\frac{v(v+1)}{2}=2(t+1)(4t+5)$ and needs to be partitioned into $2(t+1)$ positive parts. Each part represents the sum of two elements in one defining set. A result from~\cite{milenkovic2004probabilistic} shows that the number of partitions of an integer with additional constraints on the number and size of the parts can be determined asymptotically as follows. Let $V_i$ denote the number of parts of size $i$ in a randomly chosen partition, where $i=1,2,\ldots$. Assume that $V_1, V_2, \dots$ are independent and geometrically distributed random variables with parameters $\lambda, \lambda^2, \dots$, respectively, where $\lambda=\exp \left(-\pi \sqrt{\frac{1}{6 n}}\right)$. Then, the asymptotic partition count of interest can be obtained by performing all computations under the previous model. The formula for $\lambda$ is a consequence of the well-known asymptotic formula by Hardy and Ramanujan for the number of partitions of a positive integer $n$, $p(n) \sim \frac{1}{4 \sqrt{3} n} \exp \left(\pi \sqrt{\frac{2 n}{3}}\right)$. 

Note that by the definition of balanced defining sets, the smallest possible sum of a defining set is $5$ (i.e., $1+4$), and the largest possible sum is $2v-3$ (i.e., $v-3+v$), implying that $V_i\neq 0$ and has to be even for $5\leq i\leq 2v-3$. Therefore, the overall probability of drawing valid constrained partition for an integer that is the sum $s$ equals
\begin{align*}
P_s=&\sum_{x_5,\dots,x_{2v-3}}\lambda^{2\left[5x_5+6x_6+\dots+(2v-3)x_{2v-3}\right]}\prod_{i=1}^s(1-\lambda^i)\notag \\
=&\binom{\frac{v}{4}+2v-8}{2v-8}\lambda^s\prod_{i=1}^s(1-\lambda^i);\; \lambda=\exp \left(-\pi \sqrt{\frac{1}{6 s}}\right)
\end{align*}
The last line follows from the facts that $\sum_{j=5}^{2v-3}jx_j=\frac{s}{2}$ and $\sum_{j=5}^{2v-3}x_j=\frac{v}{4}$. For the sum of two elements in each defining set, there are at most $\frac{v}{2}$ possibilities. Therefore, the final asymptotic upper bound for the number of valid partitions equals $P_s\times p(s)\times \frac{(v/2)^{v/2}}{(v/2)!}$. Details of the proof are omitted.

\section{Balanced and Swap-Robust Minimal Trades} \label{sec:balanced-robust-trades} 

Among the large number of valid balanced defining sets, we now wish to select those that are most resilient to what we call popularity-swaps. We introduce the following definitions.

\begin{definition}
We say that a data chunk labeled $i$, where $i \in [v]$, experienced a \emph{popularity change} of magnitude $p$ if its label changes to a value $i_p$,  $i_p \in [v]$, so that $|i-i_p| \leq p$. Note that the most and least popular items, labeled by $1$ and $v$, can only increase or decrease their popularity, respectively and hence the rankings are not to be viewed as cyclic. Since in the presence of popularity changes the structure of a ranking (permutation) needs to be preserved, a popularity increase has to be matched by a popularity decrease. We say that $(i,j), i,j \in [v],$ experience a popularity change of magnitude $p$ if the rankings of $i$ and $j$ swap and $|i-j|=p$. 
\end{definition}
For example, the transposition $(1,2)$ corresponds to the top-ranked element becoming the second-ranked one and vice versa. The magnitude of the popularity change equals $p=1$. As another example, for $p=1$, and $v=4$, allowed popularity swaps include $(1,2),(2,3),(3,4)$, but only collections of swaps that do not involve the same element (number) can be used simultaneously. Two ``adjacent'' swaps involving one common element are considered as one ``non-adjacent'' swap that could result in larger popularity drops by our definition. For example, allowing both swaps $(1,2)$ and $(2,3)$ simultaneously would result in $1$ dropping to $3$ and $3$ rising to $1$, which is a popularity change of magnitude $p=2$. For simplicity, we focus on the case $p=1$ and only provide a sample result for $p=2$.

\begin{example} Suppose that $t=3$ and that the defining sets are defined as in the proof of Lemma~\ref{lem:baltrade}.
Based on Theorem~\ref{thm:hwang}, the following two blocks belong to $T^{(1)}$:
\begin{align*}
B_1 = \{ 1,4,6,7,9,12,14,15 \}, B_2 = \{ 2,3,5,8,10,11,13,16 \}.
\end{align*}
Consider the swaps $(1,2),(4,5),(7,8),(9,10),(12,13),$ $(15,16)$ which change the defining sets $S$ to
\begin{align*}
&S_{1} = \{2,5\}, \;S_2 = \{1,3\}; \; S_3 = \{4,7\}, \;\;  S_4 = \{6,8\}; \\
&S_5 = \{10,13\}, \;\;  S_6 = \{9,11\}; \; S_7 = \{12,15\},  \;\;  S_8 = \{14,16\}, 
\end{align*}
and the blocks $B_1$ and $B_2$ to
\begin{align*}
B_1 = \{ 2,5,6,8,10,13,14,16 \}, B_2 = \{ 1,3,4,7,9,11,12,15 \}.
\end{align*}
The discrepancy of the sum of the two blocks equals $12$, which equals $3(t+1)$. In general, the block-discrepancy can increase to $3(t+1)$ if each pair of sets $S_{2i-1}, S_{i}$ is modified by increasing/decreasing the values of two elements from one of the sets $S \in \{S_{2i-1}, S_i\}$ and decreasing/increasing the values of two elements $\{S_{2i-1}, S_i\} \setminus S$;  this results in a contribution of $3$ to the total change in the block sum for every pair of defining sets. Assuming each of the defining sets $S$ contains exactly two elements and a similar reasoning as described above, it is clear that the block-discrepancy can increase to at most $4(t+1)$ following popularity changes of magnitude $p=1$.
\end{example}
Based on the example, it appears that different choices of balanced defining sets of minimal trades may have different tolerance to popularity swaps of magnitude $p=1$. The question of interest is to determine which choices of balanced defining sets have the smallest sums of discrepancies between companion sets, termed \emph{total set discrepancy}, under worst-case popularity changes of magnitude $p=1$. 
\begin{lemma}\label{lem:lbdis} The total set discrepancy of minimal trades with $(v=4(t+1),k=2(t+1),t)$ and defining sets of cardinality $2$, under popularity changes of magnitude $p=1$, is at least $\frac{2}{3}(t+\frac{2}{3})$.
\end{lemma}
\begin{IEEEproof} We model the set of potential popularity changes of magnitude $1$ as edges of a digraph $G=(V,E)$ with $|V|=2(t+1)$. The graph $G$ has the following properties:
\begin{enumerate}
\item Each vertex $v_i \in V$ corresponds to the defining set $S_i$;
\item There is a directed arrow from vertex $v_{i_1}$ to $v_{i_2}$, denoted by $(v_{i_1}, v_{i_2})$, if $\exists s$ such that $s \in S_{i_1}$ and $s+1 \in S_{i_2}$. We allow for parallel edges but do not allow self-loops.
\end{enumerate}
Since the trade is balanced, if any defining set $S_i$ contains a pair of consecutive elements, then the elements in its companion set cannot contain a pair of consecutive integers. Since self-loops are not allowed, it follows that the number of edges in $E$ is at least $3(t+1)-1$: From the $4(t+1)-1$ possible pairs $(s,s+1)$, there can be as many as $t+1$ such pairs that belong to the same defining set.

``Edge selection'' refers to the process where we select edges $(v_{i_1}, v_{i_2})$ that correspond to the popularity of an element in $S_{i_1}$ changed to that of an element in $S_{i_2}$. In order to ensure that the edges selected correspond to popularity changes of magnitude at most
$1$, the following needs to hold: For any vertex, if an incoming directed arc into the vertex is selected, then an outgoing arc from the same vertex cannot be selected. Similarly, if for a vertex an outgoing arc is selected, no incoming arcs can be selected for that vertex.

Let $E'$ stand for a set of arcs describing a collection of popularity changes of magnitude one. For $j \in \{1,2,\ldots,2t+2\}$, $S'_{j}$ denotes the defining sets after the popularity changes induced by $E'$. It is straightforward to see that 
because the trades are balanced, the change in discrepancy after the popularity updates equals $\sum_{j \in \{1,2,\ldots,t+1\}} \left | \sum_{i \in S'_{2j}} i - \sum_{i' \in S'_{2j-1}} i' \right|$, where
\begin{align}\label{eq:discg}
\sum_{j \in \{1,2,\ldots,t+1\}} \left | \sum_{i \in S'_{2j}} i - \sum_{i' \in S'_{2j-1}} i' \right| \leq 2| E' |.
\end{align}
We exhibit an edge set $E'$ so that~\eqref{eq:discg} holds with equality. 

If $(v_{i_1}, v_{i_2}) \in E'$, where $i_1 \in \{2j_1-1, 2j_1\}$ and $i_2 \in  \{2j_2-1, 2j_2\}$, then the following hold:
\begin{enumerate}
\item For any $v_{\ell} \in V$, $(v_{\ell}, v_{i_1})$ and $(v_{i_2}, v_{\ell})$ cannot both be included in $E'$;
\item Let $S_s, s=i_1\pm 1,$ be the companion defining set for $S_{i_1}$. For any $v_{\ell} \in V$, we cannot select $(v_{s}, v_{\ell})$;
\item  Let $S_s, s=i_2\pm 1,$ be the companion defining set for $S_{i_2}$. For any $v_{\ell} \in V$, we cannot select $(v_{\ell}, v_{s})$.
\end{enumerate}

If we select $E'$ according to (1)-(3), we will select at least $1$ out of every $9$ edges in $E$. Given that $|E| \geq 3(t+\frac{2}{3})$, the resulting change in discrepancy is at least
$
2 |E'| \geq 2 \left( \frac{|E|}{9} \right) \geq \frac{6(t+\frac{2}{3})}{9}.
$
\end{IEEEproof}
Due to space limitations, we omit a similar proof for the following sharper lower bound. 
\begin{theorem}\label{th:lbdis} Provided $19 | (t+1)$, the total set discrepancy is at least $\frac{14}{19}(t+1)$.
\end{theorem}

\begin{lemma} There exists a minimal balanced trade with $v=16$ that has total set discrepancy $\leq 6$ for any collection of popularity swaps of magnitude $p=1$.
\end{lemma}
\begin{IEEEproof}
Choose the following defining sets, pictured in Figure~\ref{fig:example}, top:
\begin{align*}
&S_1 = \{1,16\}, S_2 = \{8,9\}; \; S_3 = \{4,5\}, S_4 = \{2,7\}; \\
&S_5 = \{10,15\}, S_6 = \{12,13\}; \; S_7 = \{3,14\}, S_8 = \{6,11\}.
 \end{align*}

For any popularity swaps of $p=1$ involving elements in $S_3,S_4,S_5,S_6$, one has
\begin{align}\label{eq:discp1}
\left | \sum_{j \in S'_i} j -  \sum_{j \in S'_{i+1}} j \right | \leq 2,
\end{align}
where $i \in \{{3,5\}}$. First consider $i=3$. Note that the adjacent swap $(4,5)$ does not change the balance. A simultaneous pair of adjacent swaps $(3,4)$ and $(4,5)$ is impossible, while the pair of adjacent swaps $(3,4)$ and $(5,6)$ leaves the sum unchanged. In this case, the only allowed adjacent swaps involving elements $2,7$ would be $(1,2)$ and $(7,8)$, which when applied together leave the sum of elements in $S'_4$ unchanged. The same is true for $i=5$. Hence, only one transposition per set is possible, leading to the upper bound above. Similarly, 
\begin{align}\label{eq:discp2}
\left | \sum_{j \in S'_1} j -  \sum_{j \in S'_{2}} j \right | \leq 2.
\end{align}
This follows since the only swaps that can affect $1$ and $16$ are $(1,2)$ and $(16,15)$, which in combination leave the sum unchanged. On the other hand, the sum of the second set, $S_{2}$, can change by at most $1$, either due to the transposition $(7,8)$ or the transposition $(9,10)$. As before, one can have at most one transposition per defining set, which produces the desired upper bound. A similar line of reasoning suggest that 
\begin{align}\label{eq:discp3}
\left | \sum_{j \in S'_7} j -  \sum_{j \in S'_{8}} j \right | \leq 4,
\end{align}  
since it is possible to simultaneously increase/decrease both elements in both defining sets. Nevertheless, it is straightforward to show that if the above inequality is an equality,
the second and third pair of defining sets do not change their sums, i.e., $S'_3 = \{4,6\}, S'_4 = \{3,7\}$ and $S'_5 = \{11,14\}, S'_6 = \{12,13\}$. A similar analysis can be used to examine the remaining choices for discrepancies of the first pair of defining sets. Hence, in the worst case, the discrepancy cannot exceed $6$. 
\end{IEEEproof}
\begin{figure}[t]
\centering
\includegraphics[scale=0.2]{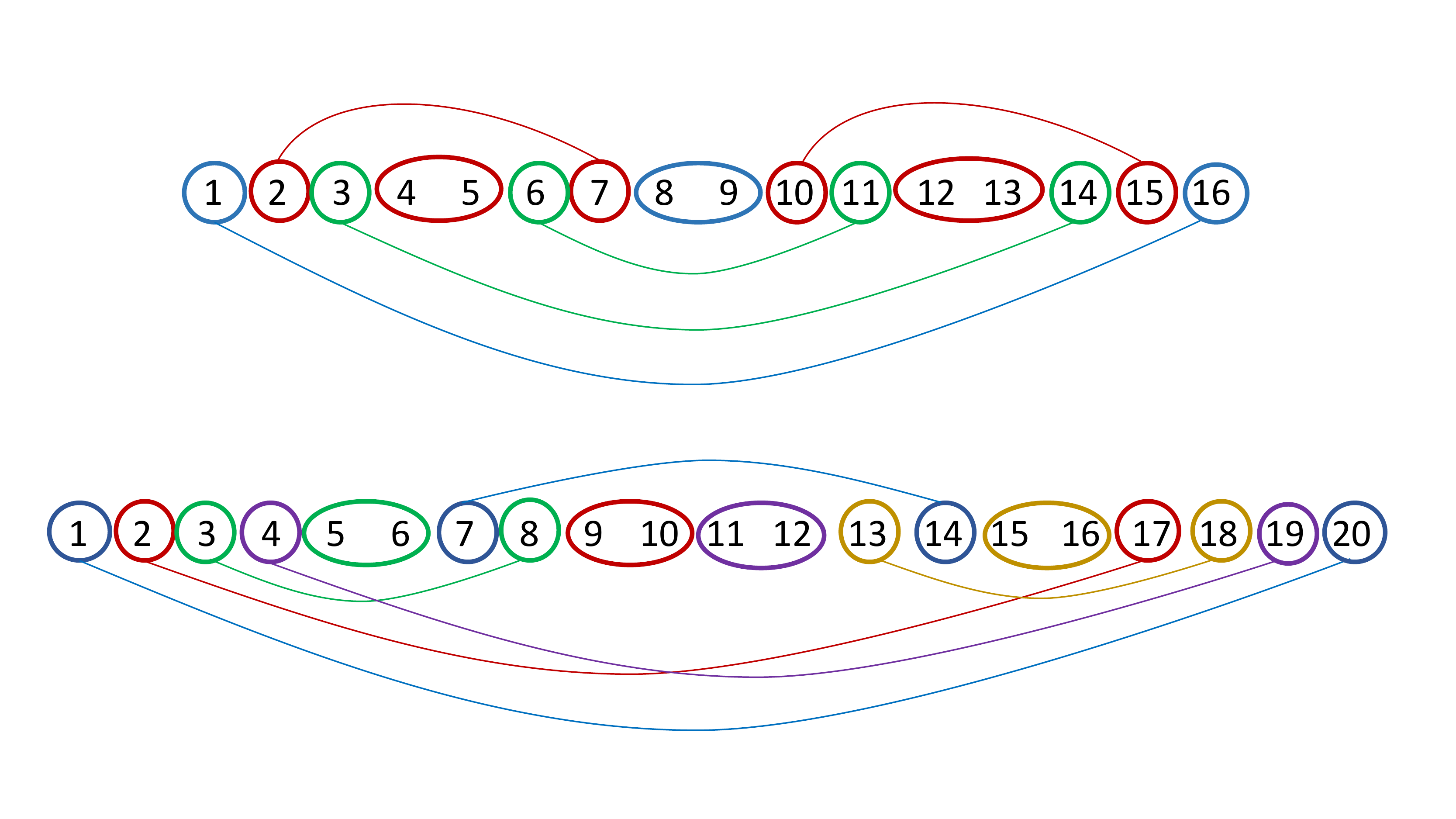}
\vspace{-0.16in}
\caption{Defining sets with smallest discrepancy for popularity changes of magnitude $p=1$ and $v=16$ (top) and $v=20$ (bottom).}
\label{fig:example}
\vspace{-0.16in}
\end{figure}
Given the large number of balanced defining sets, we only performed computer simulations for $v=12,16,20,24$ which show that the optimal total set discrepancies for these cases equal $6,6,8,10$, respectively. In particular, for $v=12$, the smallest possible change equals $2(t+1)$ as opposed to the value $3(t+1)$ described above. There are exactly $10$ collections of balanced defining sets with discrepancy $6$. Interestingly, for $v=16$ and $v=20$, the defining sets of minimum discrepancy under popularity changes of magnitude $p=1$ are unique, and depicted in Figure~\ref{fig:example}. For $v=24$, there are exactly $22$ optimal collections of defining sets.

To construct defining sets for for arbitrary values $v=4(t+1)\geq 28$, one can ``concatenate'' the structures for $v=12,16,20,24$ by subdividing the elements of the sets into largest possible groupings of $24$ or $20$ or $16$ elements and then grouping the remaining elements into one set. For example, the defining sets for $v=28$ can be constructed using the optimal patterns for defining sets with respect to elements in $\{{1,2,\ldots,16\}}$ and $\{{1,2,\ldots,12\}}$. In this case, one has to account for the discrepancy arising from swaps of elements at the ``border'' of the two groups (in this case, involving $16$ and $17$). It is easy to check that in this case, an additional discrepancy of $2$ may arise, which for the given example gives rise to a total discrepancy of $6+6+2=14=2(t+1)$. Similarly, for $v=32$, we would group the first $20$ and the last $12$ elements together and obtain a total set discrepancy $8+6+2=16=2(t+1)$. Clearly, to reduce the boundary effects, one should maximize the size of groups of integers on which the individual patterns described above are to be used. A formal statement of this result follows.

\begin{theorem}
Let $v=4(t+1)$, $m=\lfloor \frac{t+1}{6} \rfloor$ and $m'=v-24m$. Then there exists a balanced trade with block-size $k=2(t+1)$ and total set discrepancy
$$12m+\Delta=12 \lfloor \frac{t+1}{6} \rfloor+\Delta,$$
where $\Delta=-2,2,4,6,6,8$ for $m'=0,4,8,12,16,20$, respectively.
\end{theorem}
For $p=2$, we also have the following result.
\begin{theorem}
Let $v=4(t+1)$, $m=\lfloor \frac{t+1}{5} \rfloor$ and $m'=v-20m$. Then there exists a balanced trade with block-size $k=2(t+1)$ and total set discrepancy
$$26m+\Delta=26 \lfloor \frac{t+1}{5} \rfloor+\Delta,$$
where $\Delta=-4,4,8,12,18$ for $m'=0,4,8,12,16$, respectively.
\end{theorem}
The result is a consequence of the fact that for $v=12,16,20$, the optimal total set discrepancies equal $12, 18, 22$, respectively. For $v=12$, there is unique defining sets of minimum discrepancy, while there are $12,7$ optimal choices of defining sets for $v=16,20$, respectively. 

Many open problems remain, pertaining to tighter bounds on discrepancies for $p\geq 1$, constructions of (optimal) defining sets of cardinality $\geq 2$ and others.

\bibliographystyle{plain}
\bibliography{MinSumDesigns.bib}

\end{document}